\documentclass[12pt, leqno]{amsart}

\usepackage{amssymb,amsfonts,amsmath}
\usepackage{amsthm}
\usepackage{graphicx,epsfig,color}
\usepackage{enumitem}
\def\Re{{\rm{Re}}\,}
\def\Im{{\rm{Im}}\,}
\def\d{\,{\rm{d}}}
\newtheorem{theorem}{Theorem}
\newtheorem*{theorem*}{Theorem}
\newtheorem{thmx}{Theorem}

\newtheorem{proposition}{Proposition}
\newtheorem{lemma}{Lemma}

\usepackage[a4paper, total={6in, 8in}]{geometry}

\begin{document}

\title[Universality in Short Intervals]{Notes on Universality in Short Intervals and Exponential Shifts}

\thanks{\bf Dedicated to Antanas Laurin\v cikas on the occasion of his 75th birthday}

\author[Johan Andersson et al.]{Johan Andersson, Ram\= unas Garunk\v stis, Roma Ka\v cinskait\.e,\\ Keita Nakai, \L ukasz Pa\'nkowski, Athanasios Sourmelidis,\\ Rasa \& J\"orn Steuding, Saeree Wananiyakul}

\address{Johan Andersson\\ Division of Mathematics \\ Institute of Science and Technology\\ Örebro University}
\email{johan.andersson@oru.se}

\address{Ram\= unas Garunk\v stis\\ Institute of Mathematics\\
Faculty of Mathematics and Informatics\\
Vilnius University\\
Naugarduko 24, 03225 Vilnius, Lithuania}
\email{ramunas.garunkstis@mif.vu.lt}

\address{Roma Ka\v cinskait\.e\\
Institute of Mathematics\\
Faculty of Mathematics and Informatics\\ Vilnius University, Naugarduko str. 24, 03225 Vilnius, Lithuania} \email{roma.kacinskaite@mif.vu.lt}

\address{Keita Nakai\\
Graduate school of Mathematics\\ Nagoya University\\
Chikusa-Ku, Nagoya 464-8602, Japan}
\email{m21029d@math.nagoya-u.ac.jp}

\address{\L ukasz Pa\'nkowski\\ 
Faculty of Mathematics and Computer Science\\
Adam Mickiewicz University\\
Uniwersytetu Pozna\'nskiego 4, 61-614 Pozna\'n, Poland}
\email{lpan@amu.edu.pl}

\address{Athanasios Sourmelidis\\
Institute of Analysis and Number Theory\\ Graz University of technology, Steyrergasse 30, 8010 Graz, Austria}
\email{sourmelidis@math.tugraz.at}

\address{Rasa \& J\"orn Steuding\\
Department of Mathematics\\ W\"urzburg University, Emil Fischer-Str. 40, 97\,074 W\"urzburg, Germany} \email{joern.steuding@uni-wuerzburg.de}

\address{Saeree Wananiyakul\\
Department of Mathematics and Computer Science, Faculty of Science, Chulalongkorn University, 10330 Bangkok, Thailand}
\email{s.wananiyakul@hotmail.com}

\maketitle

\begin{abstract}
We improve a recent universality theorem for the Riemann zeta-function in short intervals due to Antanas Laurin\v cikas with respect to the length of these intervals. Moreover, we prove that the shifts can even have exponential growth. This research was initiated by two questions proposed by Laurin\v cikas in a problem session of a recent workshop on universality.
\end{abstract}

\section{Introduction}

In 2019, Antanas Laurin\v cikas \cite{laurincikas} proved that Voronin's celebrated universality theorem holds for shifts restricted to short intervals. More precisely, the real shifts $\tau$ for which the Riemann zeta-function $\zeta(s+i\tau)$ approximates an admissible target function (from a large class of functions) can be found in any short interval $[T,T+T^{{1/ 3}+\epsilon}]$ (for sufficiently large $T$), and the set of these shifts has positive lower density as $T\to \infty$; see Theorem \ref{Original} below. Note that the case of {\it weighted} universality is related but different; the concept of weighted universality was also introduced by Laurin\v cikas \cite{laurincikas0} in 1995.

Before we recall the precise statement we introduce some language. We say that the function ${\mathcal F}$ is {\it universal in an interval} $[T,T+H]$ if for every $\varepsilon>0$, and every admissible target function $g$, defined on some admissible set ${\mathcal K}$, the inequality
$$
\liminf_{T\to\infty}\frac{1}{ H}{\rm{meas}}\,\left\{\tau\in[T,T+H]\,:\,\max_{s\in{\mathcal K}}\vert {\mathcal F}(s+i\tau)-g(s)\vert <\varepsilon\right\}>0
$$  
holds; if here $H$ can be significantly smaller than $T$, then ${\mathcal F}$ is said to be {\it universal in short intervals}. We may restrict ourselves on the case of universality for the Riemann zeta-function ${\mathcal F}=\zeta$; Dirichlet $L$-functions can be treated analogously, and for further universal zeta- or $L$-functions the setting can be adjusted. The original universality theorem for $\zeta$ was proved by Sergei Voronin \cite{voronin} for intervals $[0,T]$ or, equivalently,  $[T,2T]$ in place of $[T,T+H]$; generalized and extended versions of Voronin's theorem were established by Steven Gonek, Bhaskar Bagchi, Axel Reich, Laurin\v cikas and others (see for example \cite{matsumoto}); the case of universality in short intervals, that is $[T,T+H]$ with $H=o(T)$, however, is new and was introduced by Laurin\v cikas \cite{laurincikas}. 
This concept was also studied for discrete universality with respect to arithmetic progressions and some other sequences \cite{lauri2}; our reasoning applies to this setting too. 
Moreover, Johan Andersson \cite{anders} has recently shown, among other things, that continuous universality in short intervals is equivalent to its discrete analogue in arithmetic progressions.

For the case of $\zeta(s)$ note that the set ${\mathcal K}$ is admissible if it is a compact subset of the strip ${\frac{1}{2}}<\sigma<1$ with connected complement, where, here and in the sequel, we write $s=\sigma+it$. Moreover, in this case a function $g\,:\,{\mathcal K}\to\mathbb{C}$ is admissible if $g$ is continuous, analytic and non-vanishing in the interior of ${\mathcal K}$. 
For the logarithm of $\zeta$, however, the target function does not need to be without zero; the non-vanishing restriction is related to the Riemann Hypothesis (see \cite[Section 9]{matsumoto}).  

In this setting Laurin\v cikas \cite{laurincikas} obtained the following result:

\begin{thmx}\label{Original}
The Riemann zeta-function is universal in short intervals $[T,T+H]$ for every $H$ satisfying
\[
T^{\frac{1}{3}}(\log T)^{\frac{26}{ 15}}\leq H\leq T.
\] 
\end{thmx}

Very recently\footnote{in a video address during a workshop on {\it Zeta-Functions, Universality, and Chaotic Operators} at CIRM, Luminy, August 7-11, 2023}, Laurin\v cikas suggested to look for improvements of this result. In particular, he asked for a proof of the above result with $H=T^\epsilon$, where $\epsilon>0$ is a fixed constant that can be arbitrarily small.
Moreover, he asked whether
\begin{align}\label{exponential}
\liminf_{T\to\infty}\frac{1}{T}\mathrm{meas}\left\{\tau\in[T,2T]:\max_{s\in\mathcal{K}}|\mathcal{F}(s+ie^\tau)-f(s)|<\varepsilon\right\}>0,
\end{align}
where $\mathcal{K}$ and $f$ are admissible and $\mathcal{F}$ is a universal function. This question probably originates from \cite{pankwoski} where a joint universality result for Dirichlet $L$-functions is proved and the vertical shifts are of the form $\tau^\alpha(\log\tau)^\beta$ for a wide variety of values for $\alpha$ and $\beta$. This result has been generalized by Laurin\v cikas {\it et al.} \cite{laurincikasetal} to vertical shifts generated by a function $\varphi(\tau)$ which have a polynomial-like behavior. 

In the subsequent sections of this note, we shall i) improve the exponent ${\frac{1}{3}}$ in Theorem \ref{Original}, ii) give an affirmative answer to the question whether $H=T^\epsilon$ is possible subject to a certain restriction of the range of universality, resp. the yet unsolved Lindel\"of Hypothesis, iii) go beyond the latter result under assumption of the open Riemann Hypothesis, and iv) give an affirmative answer to the second proposed question by introducing an alternative approach and addressing a more general problem.

We use the Landau notation $f(x) = O(g(x))$ and the Vinogradov notation $f(x) \ll g(x)$ to mean that there exists some constant $C>0$ such that $|{f(x)}| \leq C g(x)$ holds for all admissible values of $x$ (where the meaning of `admissible' will be clear from the context).

\section{Unconditional Results I --- Bourgain \& Watt}

As can be seen from \cite{laurincikas}, the one and only crucial point in proving universality for short intervals $[T,T+H]$ is a bounded mean-square for $\zeta(s)$, that is
\begin{equation}\label{ms}
\frac{1}{H}\int_T^{T+H}\vert \zeta(\sigma+it)\vert^2\d t \ll_{\sigma} 1,\quad \sigma>\frac{1}{2}.
\end{equation}
It is worth mentioning that if we can prove \eqref{ms} for some $\sigma_0>\frac{1}{2}$, then the implicit constant appearing therein is uniform on $\sigma\geq\sigma_0$ and any sufficiently large $T=T(\sigma_0)\geq0$. 
We describe this rigorously in Lemma \ref{lemma2}.
Therefore, our task to prove universality of $\zeta(s)$ in short intervals $[T,T+H]$ reduces in proving \eqref{ms} for any fixed  $\sigma>\frac{1}{2}$.
To that end the method of exponent pairs is a rather useful tool (which originates from old work by Johannes van der Corput, starting with \cite{vdc}) and we refer to the next section for a more detailed discussion.
In his paper on universality in short intervals \cite{laurincikas}, Laurin\v cikas used the following result of Aleksandar Ivi\'c \cite[Theorem 7.1]{ivic}. 

\begin{proposition}\label{aivic}
Let $(\kappa,\lambda)$ be an exponent pair satisfying $1+\lambda-\kappa\geq 2\sigma$ for $\sigma\in({1\over 2},1)$. Then (\ref{ms}) holds uniformly for
$$
T^{\kappa+\lambda+1-2\sigma\over 2\kappa+2}(\log T)^{\kappa+2\over \kappa+1}\leq H\leq T.
$$
\end{proposition}

\noindent Using the exponent pair $({9\over 26},{7\over 13})$ due to Roger Heath-Brown \cite{heath}, it follows that estimate (\ref{ms}) holds with $H=T^{23\over 70}(\log T)^{61\over 35}$. Incorporating this in \cite{laurincikas} yields a slight improvement of Theorem \ref{Original}; note that ${23\over 70}=0.32857\,142\ldots<{1\over 3}$. 

However, we can do better by making use of a recent result of Jean Bourgain and Nigel Watt \cite{watt}, namely
\begin{equation}\label{bw}
{1\over 2U}\int_{T-U}^{T+U} \vert\zeta({\textstyle{1\over 2}}+it)\vert^2\d t \ll \log T \qquad\mbox{for}\quad U=T^{{1273\over 4053}+\epsilon}.
\end{equation}
We shall show that this mean-square estimate on the critical line implies the desired mean-square bound (\ref{ms}) for the same range to the right of the critical line (up to a negligible factor $T^\epsilon$). This leads to 

\begin{theorem}\label{1}
The Riemann zeta-function is universal in short intervals $[T,T+H]$ for every $H$ satisfying
$$
T^{{1273\over 4053}}\leq H\leq T.
$$ 
\end{theorem}

\noindent Note that ${1273\over 4053}=0.31408\,832\ldots$
\medskip

The proof relies on the following lemmata.

\begin{lemma}\label{lemma1} 
Let $T\geq1$, $1\leq H \leq T$ and ${1\over 2}\leq \sigma_0<\sigma \leq 1$. 
Then
\begin{align*}\int_T^{T+H} \left| \zeta(\sigma+it) \right|^2 \d t \ll H+ \frac{ H^{2(\sigma_0-\sigma)}}{(\sigma-\sigma_0)^2 } \int_{T-\log T}^{T+H+\log T} \left| \zeta(\sigma_0+it) \right|^2 \d t.
\end{align*}
\end{lemma}
Lemma \ref{lemma1} resembles the result of Ivi\'c \cite[Lemma 7.1]{ivic} and the proofs of these results follow along the same lines.
An immediate consequence of Lemma \ref{lemma1} is that the mean-square of $\zeta(s)$ in short intervals on a vertical line implies a better mean-square estimate result to the right of this line.

\begin{lemma}\label{lemma2}
Suppose that there exists some $\theta\in(0,1)$ and $\sigma_0 \geq {1\over 2}$ such that, for any  $\epsilon>0$, $T\geq1$ and $T^{\theta+\epsilon} \leq H \leq T$, we have that
$$
\int_{T}^{T+H} \left| \zeta \left(\sigma_0+it \right)\right|^2  \d t \ll_\epsilon H T^{\epsilon}.
$$
Then, for any $\eta\in(0,\frac{1}{2})$ and any function $\omega(T) \to 0^+$ as $T \to \infty$, we have that
$$
\int_{T}^{T+H}  \left| \zeta \left(\sigma+it \right)\right|^2  \d t \ll_{\eta,\omega} H,
$$
valid for any $T>0$, $T^{\theta-\omega(T)} \leq H \leq T$ and $\sigma_0+\eta\leq\sigma\leq1$.
\end{lemma}

\noindent {\bf Proof of Lemma \ref{lemma2}.} 
It is sufficient to prove this for $H=T^{\theta-\omega(T)}$.
To that end, let $\epsilon=\min(1-\theta,\eta\theta)$.
By Lemma \ref{lemma1} we have that
\begin{align*}  
 \int_T^{T+H} \left| \zeta(\sigma+it) \right|^2 \d t &\ll_\eta H+ H^{2(\sigma_0-\sigma)}
\int_{T-\log T}^{T+H+\log T} \left| \zeta(\sigma_0+it) \right|^2 \d t \\ 
&\ll_\eta H+  H^{-2\eta} \int_{T-\log T}^{T-\log T+T^{\theta+\epsilon}} \left| \zeta(\sigma_0+it) \right|^2 \d t.
\end{align*}
In view of our assumptions, it follows that the right hand side of the above relation is bounded by
\begin{align*}
 \ll_\eta H + H^{-2\eta} T^{\theta+ \epsilon}\ll_\eta H+T^{-2\eta(\theta-\omega(T))+\theta+\eta\theta}\ll_\eta H\biggl(1+T^{-\eta\theta+(1+2\eta)\omega(T)}\biggr).
\end{align*}
Since $\omega(T) \to 0^+$ as $T \to \infty$, we deduce that the right hand side of the above relation is $O_{\eta,\omega}(H)$ and the lemma follows.\qed
\medskip

We return now to the 

\noindent{\bf Proof of Lemma \ref{lemma1}.} 
%Let $\frac 1 2 \leq \sigma_0 <\sigma < \frac 3 2$ and $t \geq 3$.
Let 
\begin{gather*} 
 \zeta_H(s)=\sum_{n=1}^\infty n^{-s} e^{-n/H}
\end{gather*}
be a smooth truncation of the Riemann zeta-function's Dirichlet series. By a variant of Perron's formula (see for example \cite[(4.60)]{ivic} we have that
\begin{gather*}
\zeta_H(s) =\frac 1 {2 \pi i} \int_{2-\infty i}^{2+\infty i} \Gamma(z) \zeta(s+z) H^z \d z.
\end{gather*}
By moving the integration path to the line $-\frac{1}{2}<\Re(z)=\sigma_0-\sigma<0$ we pick up residues from the Gamma-function at $z=0$ and from the zeta-function at $z=1-s$ and obtain %$ and pick up $\zeta(s)$ from the residue of the Gamma-function at   $z=0$ and $\Gamma(1-s)H^{1-s}$ from the residue of the zeta-function at $s+z=1$. 
\begin{gather} \label{a2}
\zeta_H(s) =\frac 1 {2 \pi i} \int_{\sigma_0-\sigma-\infty i}^{\sigma_0-\sigma+\infty i} \Gamma(z) \zeta(s+z) H^z \d z + \zeta(s)+\Gamma(1-s)H^{1-s}.
\end{gather}
We now use the change of variables $\tau= -i(z+\sigma-\sigma_0)$  and divide the integral in \eqref{a2} into one part with $|\tau| \leq \log T$ and one part with $|\tau| \geq \log T$. After rearranging the terms,   and multiplying out a factor $H^{\sigma_0-\sigma}$ from the integrals we obtain
\begin{gather}\label{zetadef}
  \zeta(s)=\zeta_H(s)-H^{\sigma_0-\sigma} R_H(s)- H^{\sigma_0-\sigma} E_H(s),%\Gamma(1-s)H^{1-s}
  \end{gather}
{where}
\begin{gather*}
 R_H(s)= \frac {1} {2 \pi} \int_{-\log T}^{\log T}  \Gamma(\sigma_0-\sigma+i\tau)
\zeta(\sigma_0+it+i \tau) H^{i \tau} \d\tau, \\ \intertext{and}
E_H(s)= \Gamma(1-\sigma-it)H^{1-\sigma_0+it}+ \frac{1}{2 \pi} \int_{|\tau| \geq \log T} \Gamma(\sigma_0-\sigma+i\tau)
\zeta(\sigma_0+it+i \tau) H^{i \tau} \d \tau.
\end{gather*}
In view of Stirling's formula \cite[(A.34)]{ivic}
\begin{gather}\label{Stirl}
\left| \Gamma(x+it)\right| \ll e^{-|t|}(x+|t|)^{-1}, \quad 0<x\leq1,\,t\in\mathbb{R},
\end{gather}
and the bound
\begin{gather}
 \zeta(\sigma_0+it) \ll |t|^{\frac{1}{4}}, \quad |t| \geq 1,
\end{gather}
it follows that
\begin{align} \label{irt}
\begin{split}
  E_H(s)&\ll e^{-|t|}H^{1-\sigma_0}+\int_{\log T}^{\infty}e^{-\tau}(|t|+\tau)^{1/4}\mathrm{d}\tau\\
  &\ll e^{-|t|}H^{1-\sigma_0}+e^{-\pi|t|/4}|t|^{1/4},\qquad\qquad\qquad|t|\geq1.
  \end{split}
\end{align}
Therefore, in view of \eqref{zetadef}, we obtain that
\begin{gather} \label{eq19}
 \int_T^{T+H} |\zeta(\sigma+it)|^2 \d t\ll I_1+H^{2 (\sigma_0-\sigma)} I_2+o(1), \end{gather}
 where
 \begin{gather*}
I_1=\int_T^{T+H} |\zeta_H(\sigma+it)|^2 \d t\quad\text{ and }\quad I_2 =\int_T^{T+H} |R_H(\sigma+it)|^2 \d t.
\end{gather*}
 By the Montgomery-Vaughan inequality \cite[Corollary 3]{MontgomeryVaughan}   we get that 
\begin{gather} \label{montvau}
  I_1 \ll H\sum_{n=1}^\infty n^{-2 \sigma} e^{-2n/H} +\sum_{n=1}^\infty {n^{1-2 \sigma}} e^{-2n/H}\ll H,\quad H\geq1.
\end{gather}

We now consider the integral $I_2$.  By $|R_H(s)|^2 =R_H(s) \overline{R_H(s)}$ we can write $I_2$ as
\begin{multline} \label{alberg}
  I_2=\frac 1 {(2 \pi )^2} \int_{T}^{T+H} \int_{-\log T }^{\log T} \int_{-\log T }^{\log T }  H^{i (\tau_1-\tau_2)}  \Gamma(\sigma_0-\sigma +i \tau_1) \times \\ \times \Gamma(\sigma_0-\sigma- i \tau_2) \zeta(\sigma_0+i(t+\tau_1)) \zeta(\sigma_0-i(t+\tau_2)) \d \tau_1 \d \tau_2 \d t
\end{multline}
We change the integration order so that $t$ is the innermost variable of integration. The integral in $t$ may now be estimated by the Cauchy-Schwarz inequality
\begin{align*}
&\left|\int_T^{T+H}    \zeta(\sigma_0+i(t+\tau_1)) \zeta(\sigma_0-i(t+\tau_2)) \d t\right|^2 \\
& \qquad \qquad \leq  \int_{T+\tau_1}^{T+H+\tau_1} \left|  \zeta(\sigma_0+it)\right|^2\d t   \int_{T+\tau_2}^{T+H+\tau_2} \left|  \zeta(\sigma_0+it)\right|^2\d t \\ & \qquad \qquad \leq  \left(\int_{T-\log T}^{T+H+\log T} \left|  \zeta(\sigma_0+it)\right|^2\d t \right)^2,
\end{align*}
where the last inequality follows from $-\log T \leq \tau_1,\tau_2 \leq \log T$. 
By estimating the rest of the factors in \eqref{alberg} by their absolute values and using that the innermost integral is now independent of $\tau_1,\tau_2$ we get that
\begin{gather} \label{rty}
 I_2 \ll\left( \frac 1 {2 \pi} \int_{-\log T}^{\log T} \left|   \Gamma(\sigma_0-\sigma+i \tau)  \right| \d\tau \right)^2   \int_{T-\log T}^{T+H+\log T} \left|  \zeta(\sigma_0+it) \right|^2 \d t
\end{gather}
and the first factor in the above relation is bounded by $(\sigma-\sigma_0)^{-2}$ in view of Stirling's formula \eqref{Stirl}.
Now the conclusion of Lemma \ref{lemma1} follows from  \eqref{eq19}, \eqref{montvau} and \eqref{rty}.\qed
\medskip

Applying Lemma 1 and Lemma 2 in combination with \eqref{bw} implies the statement of Theorem \ref{1}.  
The method of proof also shows that any mean-square bound on the critical line of the form 
\begin{equation}\label{crit}
{1\over H}\int_{T}^{T+H} \vert\zeta({\textstyle{1\over 2}}+it)\vert^2\d t \ll T^\epsilon \qquad\mbox{for}\quad H=T^\eta
\end{equation}
with some $\eta>0$ implies universality for $\zeta$ in short intervals $[T,T+H]$ for every $H$ satisfying $T^{\eta}\leq H\leq T$.
We will return to this observation in Section 4.

\section{Unconditional Results II --- Restricted Universality}

There is a variation in using Ivi\'c's Proposition \ref{aivic} which improves the exponents for a certain prize.
For this purpose we introduce the concept of restricted universality as follows: we say that the function ${\mathcal F}$ is $\sigma_0$-{\it restricted universal in an interval} $[T,T+H]$ if it is universal for the same interval and for every admissible set ${\mathcal K}$ located in the restricted strip $\sigma_0<\sigma<1$, where $\sigma_0\in({1\over 2},1)$. 
Examples of $L$-functions which are $\sigma_0$-restricted universal in $[T,2T]$ are elements of the Selberg class with large degree. Again we shall consider only the case of the Riemann zeta-function since our methods can be adjusted to more general $L$-functions as well.

As we have discussed in the previous section, in order to show that $\zeta(s)$ is $\sigma_0$-restricted universal in $[T,T+H]$, it suffices to obtain \eqref{ms} for $\sigma_0$.
If $\sigma_0$ is ``far'' from $\frac{1}{2}$ (the prize we have to pay) then we can prove \eqref{ms} for even shorter intervals.

For a better understanding we first recall some theory of exponent pairs (see \cite[Section 2.3]{ivic}). 
A pair of non-negative real numbers $(\kappa,\lambda)$ is said to be an {\it exponent pair} if $0\leq\kappa\leq\frac{1}{2}\leq\lambda\leq1$ and
$$
\sum_{B<n\leq B+h}\exp(2\pi if(n))\ll A^\kappa B^\lambda,
$$
where $A>\frac{1}{2}$, $B\geq 1, 1<h\leq B$, and $f$ is a differentiable function satisfying  
$$
f'(x)\asymp A\qquad \mbox{for}\quad B\leq x\leq 2B.
$$
For some {\it advanced} exponent pairs one may have to assume similar bounds for higher derivatives but this is not relevant for our application. More essential in our context is the question how to find exponent pairs.

It is not difficult to verify that $(0,1)$ and $({1\over 2},{1\over 2})$ are exponent pairs. 
In addition, if $(\kappa, \lambda)$ is an exponent pair, then the so-called {\it $A$-} and {\it $B$-processes} (or else Weyl's differencing method, resp. van der Corput's method) produce further exponent pairs, namely
$$
\left({\kappa\over 2\kappa+2},{1\over 2}+{\lambda\over 2\kappa+2}\right)\quad\mbox{and}\quad \left(\lambda-{1\over 2},\kappa+{1\over 2}\right),
$$
respectively. Iteration of these processes leads to an infinitude of new exponent pairs. For the exponent in the restricted universality due to Laurin\v cikas \cite{laurincikas} the exponent pair $(\kappa,\lambda)=({4\over 11},{6\over 11})$ has been used in Ivi\'c's Proposition \ref{aivic} with $\sigma={1\over 2}+\epsilon$. 
Note that $({1\over 2},{1\over 2})$ is not applicable since it does not meet the condition $1+\lambda-\kappa\geq 2\sigma$. Actually, $({4\over 11},{6\over 11})$ does neither result from the $A$- nor the $B$-process but from convexity. If $(\kappa_1,\lambda_1)$ and $(\kappa_2,\lambda_2)$ are exponent pairs, then also
\begin{equation}\label{conv}
(t\kappa_1+(1-t)\kappa_2,t\lambda_1+(1-t)\lambda_2)
\end{equation}
is an exponent pair for any $t\in [0,1]$. We observe that then $({4\over 11},{6\over 11})$ arises from applying $A,B$ to $({1\over 2},{1\over 2})$ to get $({2\over 7},{4\over 7})$ in combination with $({1\over 2},{1\over 2})$ and the parameter $t={12\over 33}$. Note that all three exponent pairs here give the exponent ${1\over 3}$ at $T$ but the latter exponent pair $({4\over 11},{6\over 11})$, chosen by Laurin\v cikas \cite{laurincikas}, gives a smaller exponent for the $\log$-term than $({2\over 7},{4\over 7})$.  

There are further exponent pairs known which do not arise by one of the processes above; for example the pair $({13\over 84}+\epsilon,{55\over 84}+\epsilon)$ found by Jean Bourgain \cite{bourgain} which led him to the so far best bound for the Riemann zeta-function on the critical line. In our case, however, this exponent pair does not produce a $T$-exponent below ${1\over 3}$. This is another example showing that one exponent pair might be good for one application but less good for another. 

Recently, Timothy Trudgian and Andrew Yang \cite{ty} came up with an update of Robert Rankin's approach for finding the {\it best possible} exponent pair to a given problem \cite{rankin}. We found the exponent pair $({9\over 26},{7\over 13})$ by checking the boundary of the convex set of all known exponent pairs considered in their paper (see Figure 1).

\begin{figure}[h]
\includegraphics[height=6.7cm]{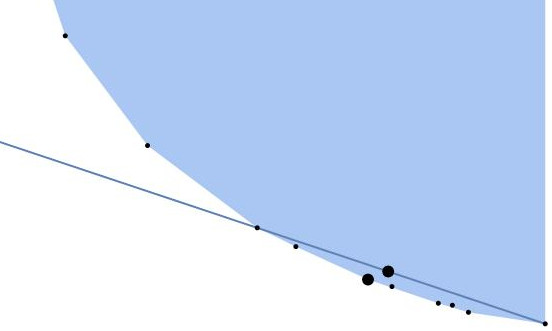}
%\caption{}
\end{figure}

\noindent The gray area is the set ${\mathcal E}$ of all known exponent pairs; the straight line passing through consists of the set of exponent pairs which yield the $T$-exponent $\frac{1}{3}$ including Laurin\v cikas' exponent pair $(\kappa,\lambda)=({4\over 11},{6\over 11})$ represented as thick dot. All points in the gray-coloured set below the line yield a $T$-exponent $<\frac{1}{3}$. The best choice, however, follows from the exponent pair $(\kappa,\lambda)=({9\over 26},{7\over 13})$ which is represented by another thick dot a little below this line. To see that one observes, by calculating the directional derivatives of $(\kappa,\lambda)\mapsto {\kappa+\lambda\over 2\kappa+2}$, that the minimum for the $T$-exponent is taken on the lower boundary of ${\mathcal E}$; since this part of the boundary consists of line segments, verifying that given ${\mathcal E}$ our choice is optimal is only a matter of a straightforward computation.  
\smallskip

We consider the condition $1+\lambda-\kappa\geq 2\sigma$ in Ivi\'c's Proposition \ref{aivic} which, for the exponent pair $({9\over 26},{7\over 13})$ implies as necessary inequality
$$
\sigma\leq {1\over 2}\left(1+{7\over 13}-{9\over 26}\right)={31\over 52}.
$$
In combination with Lemma \ref{lemma1}, this gives the existence of the mean-square (\ref{ms}) in the half-plane $\sigma>{31\over 52}$ for short intervals $[T,T+H]$ with $H\geq T^{9\over 35}(\log T)^{61\over 35}$. 
We arrive, therefore, at the following  

\begin{theorem}\label{unc}
The Riemann zeta-function is ${31\over 52}$-restricted universal in short intervals $[T,T+H]$ for every $H$ satisfying
$$
T^{9\over 35}(\log T)^{61\over 35}\leq H\leq T.
$$
Moreover, for every fixed $\epsilon\in (0,{1\over 2})$, the zeta-function is $1-\epsilon$-restricted universal in short intervals $[T,T+H]$ for every $H$ satisfying
$$
T^\epsilon\leq H\leq T.
$$
\end{theorem}

\noindent Note that ${9\over 35}=0.25714\,285\ldots$ and ${31 \over 52}=0.59615\,384\ldots$. The second statement follows from applying (\ref{conv}) to the trivial exponent pairs $(0,1)$ and $({1\over 2},{1\over 2})$ in order to derive an exponent pair $(\kappa,\lambda)$ satisfying  $0<\kappa<\epsilon$ and $1-\epsilon<\lambda<1$.

\section{Conditional Results --- Lindel\"of \& Riemann Hypotheses}

It is a folklore conjecture that every $(\epsilon,{1\over 2}+\epsilon)$ is an exponent pair. If so, then the reasoning from the previous section would imply unrestricted universality for the Riemann zeta-function in short intervals $[T,T+H]$ with $H=T^{{1\over 4}+\epsilon}$ and ${3\over 4}$-restricted universality with $H=T^\epsilon$. This conjecture also implies the Lindel\"of Hypothesis (which states that $\zeta(\sigma+it)\ll t^\epsilon$ as $t\to +\infty$ for every fixed $\sigma\geq {1\over 2}$. Assuming the latter open conjecture, however, we can show a stronger result:

\begin{theorem}\label{theo2}
If the Lindel\"of Hypothesis is true, then the Riemann zeta-function is universal in short intervals $[T,T+H]$ for every $H$ satisfying
$$
T^\epsilon\leq H\leq T.
$$ 
\end{theorem}

\noindent {\bf Proof.} 
Once more, it will suffice to prove \eqref{ms} for every fixed $\sigma>\frac{1}{2}$ and $H$ as in the theorem.
Observe that in view of Lemma \ref{lemma1}
\[
\int_T^{T+H}\vert \zeta(\sigma+it)\vert^2\d t\ll_\sigma H+H^{1-2\sigma}\int_0^{3T}|\zeta(\frac{1}{2}+it)|^2\mathrm{d}t.
\]
However, the Lindel\"of Hypothesis implies that \cite[Theorem 13.2]{titch}
\[
\int_0^{3T}|\zeta(\frac{1}{2}+it)|^2\mathrm{d}t\ll_\epsilon T^\epsilon,\quad T\geq1.
\]
Thus, for $0<\epsilon<2\sigma-1$ we obtain
\[
\int_T^{T+H}\vert \zeta(\sigma+it)\vert^2\d t\ll_\sigma H
\]
and the theorem follows.
\qed
\medskip

It is remarkable, however, that even shorter intervals are possible, at least if one is willing to assume another unproven conjecture. Assuming the Riemann Hypothesis (that is the non-vanishing of $\zeta(\sigma+it)$ for $\sigma>{1\over 2}$), Ayyadurai Sankaranarayanan \& Kotyada Srinivas \cite{sankar} proved for 
$$
{1\over 2}+{2A_1\over \log\log T}\leq \sigma\leq 1-\delta 
$$
and $\exp((\log T)^{2-2\sigma})\leq H\leq T$ with sufficiently large $T$ that
$$
{1\over H}\int_T^{T+H}\vert \zeta(\sigma+it)\vert^2\d t=\zeta(2\sigma)+O\left(\exp\Big(-A_2{(\log T)^{2-2\sigma}\over \log\log T}\Big)\right),
$$
where the implicit constant depends on $\delta>0$ and $A_2>0$ depends on $A_1>0$.
Since the lower bound decreases to $\lim_{T\to \infty}{1\over 2}+{2A_1\over \log\log T}={1\over 2}$, we obtain 

\begin{theorem}
If the Riemann Hypothesis is true, then the Riemann zeta-function is universal in short intervals $[T,T+H]$ for every $H$ satisfying
$$
\exp((\log T)^{1-\epsilon})\leq H\leq T. 
$$
\end{theorem}

\noindent Note that $\exp(((\log T)^{2-2\sigma_0})=o(T^\epsilon)$ for every $\sigma_0>{1\over 2}$.

\medskip 

It is worth mentioning that Yoonbok Lee \cite{lee} replaced the role of the mean-square in his proof of universality for Hecke $L$-functions by a density theorem or a density hypothesis. More precisely, he proved (see \cite[Theorem 3]{lee}) that the logarithm of a given Hecke $L$-function can be well approximated in the mean on the vertical segment $\sigma+it$ with $t\in[T,T+H]$ by a suitable Dirichlet polynomial, provided $T^a\leq H\leq T$, $\frac{1}{2}<a\leq1$ and
\begin{equation}\label{eq:zeroDensity}
N_L(\sigma, T, T+H) \ll H \left(\frac{H}{\sqrt{T}}\right)^{c(\frac{1}{2}-\sigma)}\log T,\qquad\text{uniformly for $\sigma\geq \frac{1}{2}$};
\end{equation}
here $N_L(\sigma, T, T+H)$ counts the number of zeros of $L(s)$ in the region $\Re(s)>\sigma$, $T<\Im(s) = t<T+H$. Then, using a straightforward modification of Voronin's proof of universality, he succeeded in showing that \eqref{eq:zeroDensity} with $H=T$ implies universality for $L(s)$ in the entire strip $\frac{1}{2}<\Re(s)<1$. By the same argument, one can show that \eqref{eq:zeroDensity} implies also universality for $L(s)$ in short intervals $[T, T+H]$. In particular, one can prove unconditionally that the Riemann zeta-function is universal in short intervals $[T,T+H]$, since \eqref{eq:zeroDensity} with $L(s)=\zeta(s)$ is known due to Atle Selberg \cite{selberg1}.  Thus, improving Selberg's zero density estimate is an alternative strategy to prove universality in short intervals.  The best result in this direction seems to be due to Ramachandran Balasubramanian \cite[Theorem 6]{balu}, who proved that for every $H$ satisfying $T^{\frac{27}{82}}\leq H\leq T$ 
\begin{equation}\label{eq:balu}
	N_\zeta(\sigma,T,T+H) \ll H\cdot  H^{\frac{2(\frac{1}{2}-\sigma)}{3-2\sigma}}(\log T)^{100},\qquad\text{uniformly for $\sigma\geq \frac{1}{2}$.}
\end{equation}
Let us note that the only reason why Theorem 3 in \cite{lee} holds only for $H\geq T^a$ with $\frac{1}{2}<a<1$ is the presence of the factor $\frac{H}{\sqrt{T}}$ in the density hypothesis \eqref{eq:zeroDensity}. Hence, one can adopt Lee's argument, essentially replacing \eqref{eq:zeroDensity} by \eqref{eq:balu}, to improve Theorem 3 from \cite{lee} and, in consequence, show universality of the Riemann zeta-function in short intervals $[T,T+H]$ with $H=T^{\frac{27}{82}+\varepsilon}$ for every $\varepsilon>0$. What is more, if one extends \eqref{eq:balu} to $H=T^\delta$ with $0<\delta<\frac{27}{82}$, then immediately we will be able to prove universality in short intervals for the same $H$.

\section{Exponential Shifts}  
We conclude this note by answering positively the second question of Laurin\v cikas regarding \eqref{exponential}.
To that end we define the class $\Phi\subseteq C^1(0,+\infty)$ of functions $\phi(\tau)>0$  with $\psi(\tau):=\phi'(\tau)/\phi(\tau)$ satisfying for any sufficiently large $T>0$  and any $\tau\geq T$ the following properties:
\begin{enumerate}[label=(\roman*)]
    \item $\phi'(\tau)$ is an increasing function with
        $\phi'(\tau+{1}/{\psi(\tau)})\ll\phi'(\tau)$,
    \item $\psi(\tau)$ is an increasing function with $A^{-1}\leq\psi(\tau)\leq\psi(\tau+1/\psi(\tau))\leq A+ \psi(\tau)$ for some absolute constant $A>0$, or a decreasing function with $\psi(\tau)\geq B/\tau$ for some absolute constant $B>0$.
\end{enumerate}
 The above definition is tailor-made for:
 \begin{enumerate}
 \item polynomials $\phi(\tau)$ of degree $\geq1$,
     \item functions $\phi(\tau)\asymp\alpha^{p(\tau)}$ for some polynomial $p(\tau)$ and $\alpha>1$,
     \item functions $\phi(\tau)\asymp\alpha^{\beta^{p(\tau)}}$ for some polynomial $p(\tau)$ and $\alpha,\beta>1$
 \end{enumerate}
 and so on. 
 More examples can be generated from the above cases by multiplying the terms of polynomials with powers of logarithms.
 
 Having introduced the class $\Phi$, we shall establish some basic properties of its elements.
 The mean-value theorem and axiom $(i)$ imply for any $C>0$ that
	\[
	\frac{\psi(\tau)}{C}\left[\phi\left(\tau+\frac{C}{\psi(\tau)}\right)-\phi(\tau)\right]\geq\phi'(\tau)=\psi(\tau)\phi(\tau)
	\]
	and, hence,
	\begin{align}\label{someinequality}
	\phi\left(\tau+\frac{C}{\psi(\tau)}\right)\geq(C+1)\phi(\tau),\quad\tau\in[T,2T].
	\end{align}

In the case when $\psi(\tau)$ is an increasing function, we set for any integer $k\geq0$ 
\[
T_0=T\quad \text{ and }\quad T_k=T_{k-1}+\frac{1}{\psi(T_{k-1})}.
\]
It then follows  by induction  from axiom $(ii)$ that $A^{-1}\leq\psi(T_k)\leq\psi(T)+kA$, $k\geq0$.
Therefore, the series $\sum_{k}1/\psi(T_{k-1})$ diverges to $+\infty$ and each of its terms contribute at most $A$.
Hence, there is $K=K(T,A)\in\mathbb{N}$ such that
\begin{align}\label{secondindequality}
	\sum_{k\leq K}\frac{1}{\psi(T_{k-1})}=T+O(1)\quad\text{ and }\quad T_K=T_0+\sum_{0\leq m\leq K-1}\frac{1}{\psi(T_m)}=2T+O(1).
\end{align} 
In conclusion, if $\psi(\tau)$ is an increasing function, then we can find points $T_0,T_1,\dots, T_K$ which form, up to an $O(1)$ error, a partition of the interval $[T,2T]$.
\begin{theorem}
Let $\mathcal{K}$ and $f$ be admissible and $\varepsilon>0$. 
If $\phi\in\Phi$, then
\[
\liminf_{T\to\infty}\frac{1}{T}\mathrm{meas}\left\{\tau\in[T,2T]:\max_{s\in\mathcal{K}}|\zeta(s+i\phi(\tau))-f(s)|<\varepsilon\right\}>0.
\]
\end{theorem}

\noindent {\bf Proof.}
 We will mainly need Voronin's universality theorem in the following form:
\begin{align}\label{Voroninimproved1}
\mathrm{meas}\left\{\tau\in [T,DT]:\max_{s\in\mathcal{K}}|\zeta(s+i\tau)-f(s)|<\varepsilon\right\}\gg T
\end{align}
for a fixed $D>1$ and any sufficiently large $T>0$.
 We will also employ a classic result from measure theory \cite[Theorem 7.26 or page 156]{rudin}:

\begin{proposition}\label{difftrans} Suppose $\phi:[a,b]\to[\alpha,\beta]$ is absolutely continuous, monotonic, $\phi(a)=\alpha$, $\phi(b)=\beta$, and $f\geq0$ is Lebesgue measurable. Then
\begin{align*}
\int_{\alpha}^\beta f(t)\d t=\int_a^bf(\phi(x))\phi'(x)\d x.
\end{align*}
\end{proposition}

We set now 
	\begin{align*}
	\mathcal{E}&:=\left\{\tau\geq0:\max_{s\in \mathcal{K}}|\zeta(s+i\tau)-f(s)|<\varepsilon\right\},\\
 S_T&:=\left\{\tau\in[T,2T]:\max_{s\in \mathcal{K}}|\zeta(s+i\phi(\tau))-f(s)|<\varepsilon\right\}
	\end{align*}
	and $1_\mathcal{E}$ to be the characteristic function of $\mathcal{E}$.
 
 imply If $\psi(\tau)$ is a decreasing function then axiom $(i)$ and Proposition \ref{difftrans} imply
 \[
\mathrm{meas}(S_T)\geq\frac{1}{\phi'(2T)}\int_{T}^{2T}1_\mathcal{E}(\phi(\tau))\phi'(\tau)\mathrm{d}\tau=\frac{1}{\phi'(2T)}\int_{\phi(T)}^{\phi(2T)}1_\mathcal{E}(\tau)\mathrm{d}\tau.
 \]
Observe that $\phi(2T)\geq\phi(T+B/\psi(T))\geq(1+B)\phi(T)$ from \eqref{someinequality}.
Therefore, by \eqref{Voroninimproved1} with $D=1+B$ we obtain that
\[
\mathrm{meas}(S_T)\gg\frac{1}{\phi'(2T)}\int_{\frac{\phi(2T)}{1+B}}^{\phi(2T)}1_\mathcal{E}(\tau)\mathrm{d}\tau\gg\frac{\phi(2T)}{\phi'(2T)}\gg\frac{1}{\psi(2T)}\gg T.
\]

If now $\psi(\tau)$ is an increasing function then we construct the partition $T_0,\dots,T_K$ of $[T,2T]$ as described above Theorem 5 and we see that
	\begin{align*}
	\mathrm{meas}(S_T)=\sum_{k\leq K}\int_{T_{k-1}}^{T_k}1_\mathcal{E}(\phi(\tau))\mathrm{d}\tau+O(1).
	\end{align*}
	Once more from  axiom $(i)$ and Proposition \ref{difftrans} it follows that
	\[
	\phi'(T_k)\int_{T_{k-1}}^{T_k}1_\mathcal{E}(\phi(\tau))\mathrm{d}\tau\geq\int_{\phi\left(T_{k-1}\right)}^{\phi(T_k)}1_\mathcal{E}(\tau)\mathrm{d}\tau.
	\]
 Observe that $\phi(T_k)\geq2\phi(T_{k-1})$ from \eqref{someinequality}.
 Therefore, by \eqref{Voroninimproved1} with $D=2$, axiom $(i)$ and relation \eqref{secondindequality}
we obtain that
    \begin{align*}
    \mathrm{meas}(S_T)\gg\sum_{k\leq {K}}\frac{1}{\phi'\left(T_{k}\right)} \int_{\phi\left(T_{k-1}\right)}^{2\phi(T_{k-1})}1_\mathcal{E}(\tau)\mathrm{d}\tau\gg\sum_{k\leq {K}}\frac{\phi\left(T_{k-1}\right)}{\phi'\left(T_{k-1}\right)}
    \gg\sum_{k\leq K}\frac{1}{\psi(T_{k-1})}
    \gg T,
    \end{align*}
    which concludes the proof of the theorem.\qed
\medskip

{\bf Acknowledgement.} R. Garunk\v stis is funded by the Research Council of Lithuania (LMTLT), agreement No. S-MIP-22-81. \L.~Pa\'nkowski is partially supported by the grant no. 2021/41/B/ST1/00241 from the National Science Centre.
A. Sourmelidis is supported by the Austrian Science Fund (FWF) project number M 3246-N.

\end{document}